\documentclass[12pt]{amsart}
\usepackage{epsf}

\newcommand{\pfpf}[1]{\medskip\par\noindent{\bf {#1}}}
\newcommand{\clcl}[1]{\medskip\par\noindent{\it {#1}:}}

\newcommand{\caca}[1]{\medskip\par\noindent{\bf {#1}:}}
\newcommand{\stst}[1]{\medskip\par\noindent{\bf {#1}:}}
\newcommand{\defglo}[1]{{\em {#1}}}
\newcommand{\mvmv}[1]{\medskip\par\noindent{\textsc {{#1}}:}}

\newcommand{\de}{\partial}
\newcommand{\sm}{\setminus}
\newcommand{\interior}{\mbox{int}}
\newcommand{\dv}{\de_v}
\newcommand{\F}{\Sigma}
\newcommand{\X}{X}

\newcommand{\Z}{Z}
\newcommand{\R}{I \! \! R}

\begin{document}

\title{A Standard Form for Incompressible Surfaces in a Handlebody}
\author{Linda Green}
\address{MSRI \\ 1000 Centennial Drive \\ Berkeley \\ CA  94720}
\email{lindag@msri.org}
\thanks{Research at MSRI is supported in part by NSF grant DMS-9022140}

\begin{abstract}
Let $\F$ be a compact surface and let $I$ be the unit interval.  This
paper gives a standard form for all 2-sided incompressible surfaces in the
3-manifold $\F \times I$.  Since $\F \times I$ is a handlebody when 
$\F$ has boundary, this standard form applies to incompressible surfaces 
in a handlebody.  
\end{abstract}

\maketitle

\section{Introduction}

Let $M$ be a $3$-dimensional manifold and let $X \subset M$ be a
properly embedded surface.  A \defglo{compression disk} for $\X
\subset M$ is an embedded disk $D \subset M$ such that $\de D \subset
\X$, $\interior(D) \subset (M \sm \X)$, and $\de D$ is an essential
loop in $\X$.  The surface $\X \subset M$ is \defglo{incompressible}
if there are no compression disks for $\X \subset M$ and no component
of $\X$ is a sphere that bounds a ball.  If $\X \subset M$ is
connected and 2-sided, then $\X$ is incompressible if and only if the
induced map $\pi_1(\X) \rightarrow \pi_1(M)$ is injective and $\X$ is
not a sphere that bounds a ball.  See, for example, \cite[Chapter
6]{Hempel}.

Let $\F$ be a compact surface and let $I$ be the unit interval $[0,
1]$.  The manifold $\Sigma \times I$ is foliated by copies of $I$,
which can be thought of as vertical flow lines.  

This paper shows that every properly embedded, 2-sided, incompressible
surface in $\F \times I$ can be isotoped to a standard form, called
``near-horizontal position''.  A surface in near-horizontal position
is transverse to the flow on $\Sigma \times I$ in the $I$ direction,
except at isolated intervals, where it coincides with flow lines.
Near each of these intervals, $\F$ looks like
a tiny piece of a helix, with the interval as its core.  A surface in
near-horizontal position can be described combinatorially by listing
its boundary curves and the number of times it crosses each line in a
certain finite collection of flow lines.  When $\Sigma$ is a compact
surface with boundary, then $\Sigma \times I$ is a handlebody. So this
paper applies, in particular, to incompressible surfaces in a
handlebody.

\section{Notation}

Throughout this paper, if $E$ is a topological space, then $|E|$
denotes the number of components of $E$.  The symbol $I$ denotes the
unit interval $[0, 1]$.  If $M$ is a manifold, then $\de M$ refers to
the boundary of $M$ and $\interior(M)$ refers to its interior.  The
symbol $\Sigma$ refers to a compact surface, and $p$ denotes the
projection map $\Sigma \times I \rightarrow \Sigma$.

The surface $\X \subset M$ is a proper embedding if $\interior(\X)
\subset \interior(M)$, $\de \X \subset \de M$, and the intersection of
$X$ with a compact subset of $M$ is a compact subset of $X$.  The map
$G: \X \times I \rightarrow M$ is a proper isotopy between $G|_{\X
\times 0}$ and $G|_{\X \times 1}$ if for all $t \in I$, $G|_{\X \times
t}$ is a proper embedding.  Unless otherwise stated, all surfaces
contained in 3-manifolds are properly embedded and all isotopies are
proper isotopies.

A connected surface $\X \subset M$ is \defglo{boundary parallel} if $\X$
separates $M$ and there is a component $K$ of $M \sm \X$ such that
$(K, \X)$ is homeomorphic to $(\X \times I, \X \times 0)$.

\section{Definition of near-horizontal position}

Let $\X$ be a surface in $\F \times I$.   Let 
$C \subset \F$ be the union of loops and arcs $p(\X \cap (\F
\times 1))$, let $C' = p(\X \cap (\F \times 0))$, and let $B
= \X \cap (\de \F \times I)$.  

$\X$ is in \defglo{near-horizontal position} if 
\begin{enumerate}
\item $C$ and $C'$ intersect transversely,  
\item $p|_B: B \rightarrow \de \F$ is injective, 
\item $p |_{\X \sm p^{-1}(C \cap C')}$ is a local homeomorphism, and
\item for any point $z \in {C \cap C'}$, there is a neighborhood $U
\subset \F$ of $z$ 
such that $p^{-1}(U) \cap \X$ either looks like the region
\par\noindent
\begin{figure}[h!tbp]
\epsfysize=1.0in
\centerline{
\epsfbox{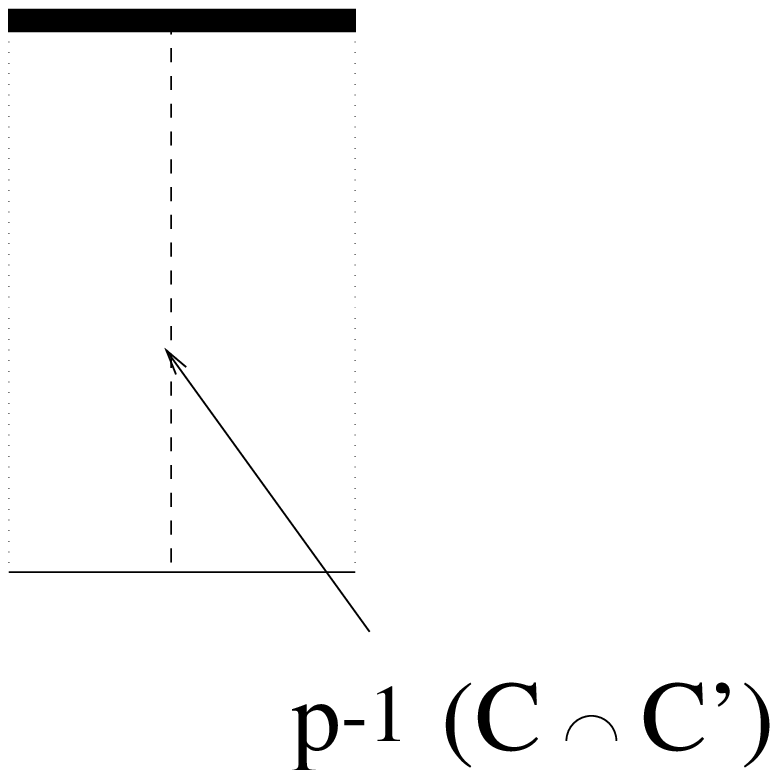}}
\end{figure}
\par\noindent
or else looks like a union of regions of the form
\begin{figure}[h!tbp]
\epsfysize=0.7in
\centerline{
\epsfbox{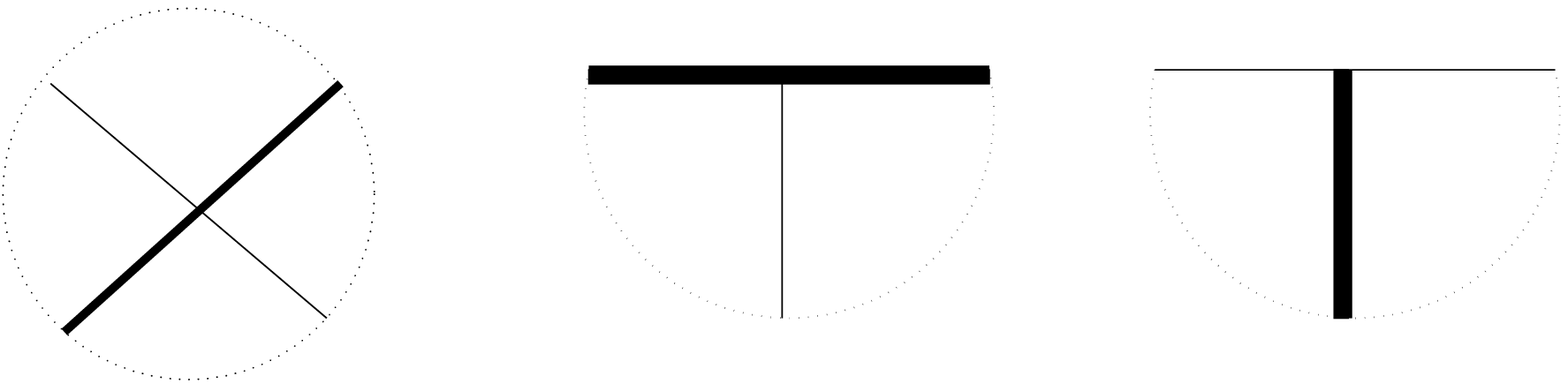}}
\end{figure}
\end{enumerate}
\par\noindent Here $p^{-1}(C) \cap \X$ is drawn with thick lines and
$p^{-1}(C') \cap \X$ is drawn with thin lines.  Pieces of the
boundaries of neighborhoods that are not part of $p^{-1}(C) \cap \X$
or $p^{-1}(C') \cap \X$ are drawn with dotted lines.  I will call the
vertical line of $p^{-1}(C \cap C')$  in the first picture a
\defglo{vertical twist line}.
\medskip

\label{section-twistsection}

\section{Examples of surfaces in near-horizontal position}

A surface $X$ in near-horizontal position determines two unions of
properly embedded arcs and loops in $\F$: $C$ and $C'$ defined above,
and a union of 
disjoint arcs of $\de \F$:  $p(B)$.  By definition of near-horizontal 
position, 
$C$, $C'$, and $p(B)$ have the following properties:
\begin{enumerate}
\item $C$ and $C'$ intersect transversely.
\item $\de (p(B)) = \de C \cup \de C'$
\end{enumerate}
Furthermore, a surface $X$  in near-horizontal position determines a 
function
\[N: \mbox{ components of $\X \sm (C \cup C')$} \longrightarrow \mbox{ non-negative integers}\]
such that for each component $r$ of $\F \sm (C \cup C')$ and each
point $y \in \interior(r)$, $N(r)$ counts the number of times $X$
intersects $y \times I$.  It is easy to check that $N$ has the
following properties:
\begin{enumerate}
\setcounter{enumi}{2}
\item If $r$ is a component of $\F \sm (C \cup C')$ with an edge in
$p(B)$, then $N(r) = 1$.
\item If $r_1$ and $r_2$ are two components which meet along an arc of $C$ 
or an arc of $C'$, then $|N(r_1) - N(r_2)| = 1$.
\item If $r_1$, $r_2$, $r_3$, and $r_4$ meet at a common vertex, then
either the set $\{N(r_1), N(r_2), N(r_3), N(r_4)\}$ contains $3$ distinct
numbers,  or else \linebreak[3] 
$\{N(r_1), N(r_2), N(r_3), N(r_4)\} = \{0, 1\}$.
\end{enumerate}

The second possibility occurs if and only if the common vertex is the
projection of a vertical twist line.

\medskip

Conversely, suppose that $C$ and $C'$ are unions of arcs and loops in
$\F$, that $p(B)$ is a union of disjoint arcs of $\de \F$, and that 
\[N: \mbox{ components of $\F \sm (C \cup C')$} \longrightarrow \mbox{ non-negative integers}\]
is a numbering scheme satisfying conditions (1) - (5) above.  Then the
information $(C, C', p(B), N)$ determines a unique surface in
near-horizontal position.

\begin{figure}[htbp]
\epsfysize=4in
\centerline{\epsfbox{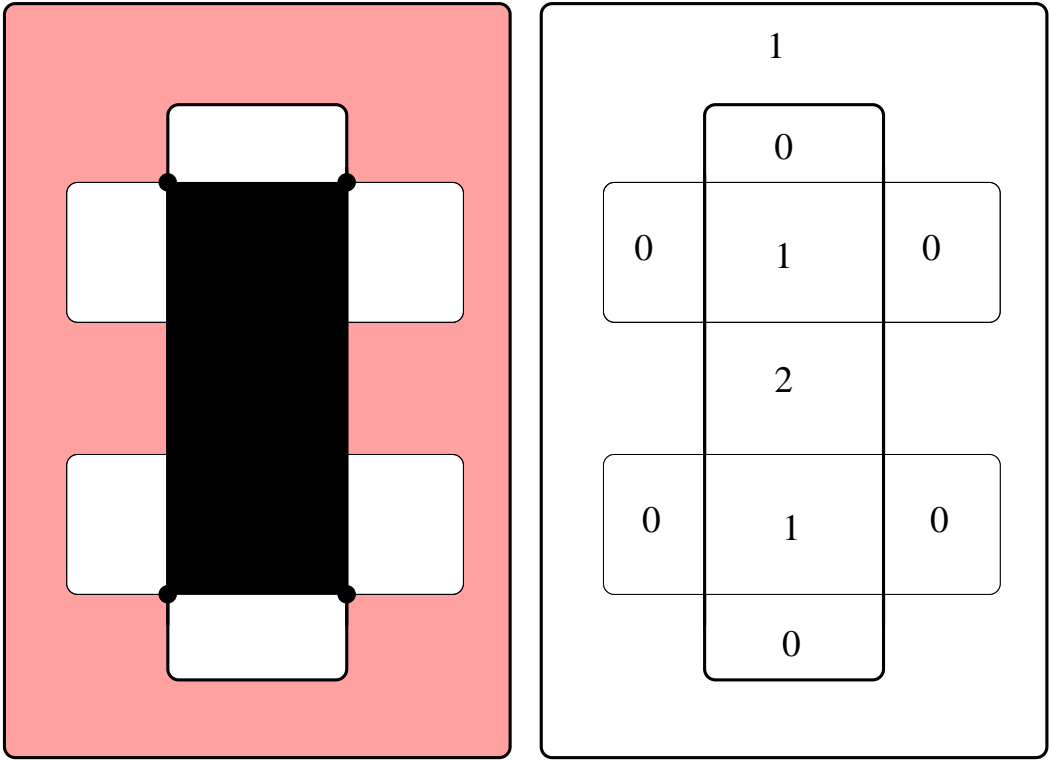}} \label{fig-example}
\caption{An example of a surface in near-horizontal position.}
\end{figure}

Figure~\ref{fig-example} depicts a genus $1$ surface with $4$ boundary
loops, in near-horizontal position in $\R^2 \times I$.  The projection
of the surface to $\R^2$ is drawn at left.  Vertical twist lines, which 
project to points, are marked with dots.  Loops of $C$ are drawn with 
thick lines, and loops of $C'$ are drawn with thin lines.  The corresponding
combinatorial description is given at right.  

\section{Putting incompressible surfaces in near-horizontal position}

\stst{Theorem}
{\it Suppose that $\X \subset \F \times I$ is a properly embedded, $2$-sided,
incompressible surface.  Then $\X$ is isotopic to a surface in
near-horizontal position.
}

\pfpf{Proof of Theorem:}
Notice that it is possible to isotope $\X$ (leaving it fixed outside a
neighborhood of $\de \X$), so that $\de \X$ satisfies requirements 1
and 2 in the definition of near-horizontal position.  Therefore, I
will assume that $C$ intersects $C'$ transversely and that $p|_B$ is
an embedding, where $C$, $C'$, and $B$ are defined as above.

Suppose that $\X$ contains components $\X_1, \X_2, \ldots \X_n$ that
are  boundary parallel disks.  Each $X_i$ can be isotoped so that
$p|_{\X_i}$ is an embedding and so that $p(\de \X_i)$ is disjoint from
$p(\de X \sm \de \X_i)$.  Therefore, if $\X \sm (\X_1 \cup \X_2 \cup
\cdots \cup \X_n)$ can be isotoped to near-horizontal position, 
so can $\X$.  So without loss of generality, I can assume that $\X$
has no components that are boundary parallel disks.  Since $\X$ is
incompressible, this assumption insures that 
all loops in $C$ and $C'$ are essential in $\F$.

All isotopies in the rest of the proof will leave $\de \X$ fixed.

Let $A$ be the union of vertical strips and annuli $C \times [0, 1]$
and let $A' = C' \times [0, 1]$.  First, I define a notion of
``pseudo-transverse'' and a measure of complexity for surfaces that
are pseudo-transverse to $A \cup A'$.  Then I describe three moves
which decrease this complexity.  In Step 1, I isotope $\X$ so that it
is pseudo-transverse to $A$ and $A'$.  In Step 2, I isotope $\X$ using
the three moves as many times as possible.  I then verify six claims
about the position of $\X$.  In Step 3, I isotope $\X$ so that the
projection map $p$ is injective when restricted to any arc or loop of
$\X \cap (A \sm A')$ and $\X \cap (A' \sm A)$.  In Step 4, I complete
the proof of the Theorem by isotoping $\X$ so that $p$ is locally
injective everywhere except at twist lines.

I will say that $\X \subset \F \times I$ is \defglo{pseudo-transverse} to $A
\cup A'$ if the following three conditions hold:
\begin{enumerate}
\item For any point $z \in (C - C')$ there is a neighborhood $U
\subset \F$ of $z$ such that $p^{-1}(U) \cap \X$ is a disjoint union
of regions of the form
\begin{figure}[h!tbp]
\epsfysize=0.7in
\centerline{
\epsfbox{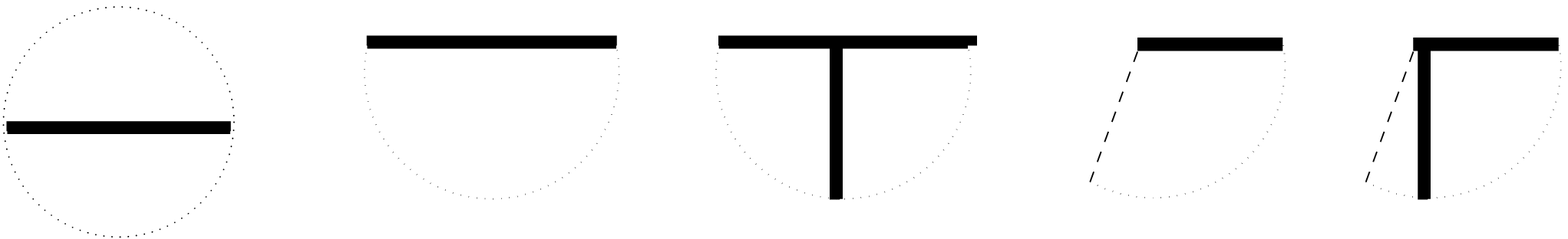}}
\end{figure}

In the last two pictures above, the dashed lines represent pieces of
$B$.  One of these two pictures occurs if and only if $z \in \de C$.
\item For any point  
$z \in (C' - C)$, there is a neighborhood U of $z$ such 
that $p^{-1}(U) \cap \X$ is a disjoint union of regions of
the form
\begin{figure}[h!tbp]
\epsfysize=0.7in
\centerline{\epsfbox{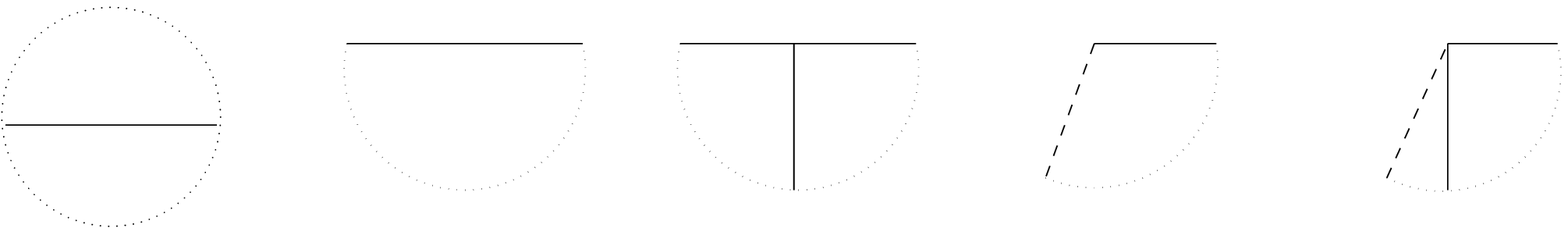}}
\end{figure}

One of the last two pictures occurs if and only if $z \in \de C'$.  
\item For any point 
$z \in (C \cap C')$, there is a neighborhood $U$ of 
$z$ such that $p^{-1}(U) \cap \X$ is either a 
single region of the form
\begin{figure}[h!tbp]
\epsfysize=1.0in
\centerline{
\epsfbox{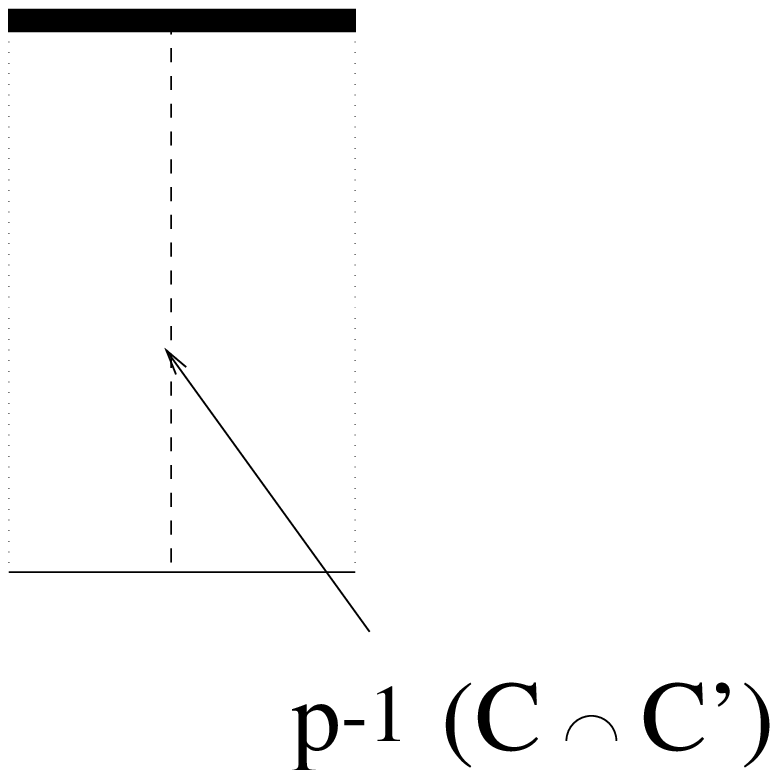}}
\end{figure}
\par\noindent
or else a disjoint union of regions of the 
form
\begin{figure}[h!tbp]
\epsfysize=0.7in
\centerline{
\epsfbox{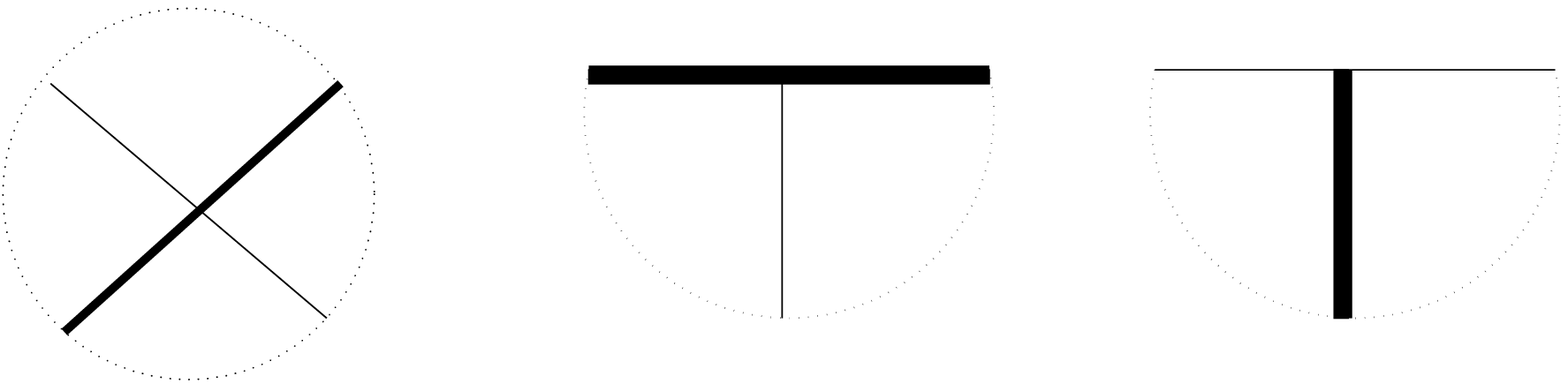}}
\end{figure}
\end{enumerate}

\medskip

Suppose $X$ is pseudo-transverse to $A \cup A'$.
Define the complexity 
of $\X$ by 
\[\zeta(\X) = (|\X \cap (A \cap A')| \mbox{, rank }
H_1(\X \cap A) + \mbox{ rank } H_1(\X \cap A')) \] ordered
lexicographically.
Consider the following three moves.

\mvmv{Move 1}  Suppose $D$ is a disk of $A$ such that $D \cap \X = 
\de D$ and $\partial D \cap A' = \emptyset$.  Then $D$ can be used to 
isotope $\X$ and decrease $\zeta(\X)$.  The isotopy leaves $\X$
in the class of pseudo-transverse embeddings.
If the roles of $A$ and $A'$ are interchanged, an analogous move is
possible.

\mvmv{Explanation of Move 1} Since $\X$ is incompressible in $\F \times I$,
$\de D$ bounds a disk $D'$ in $\X$.  The set $D \cup D'$ forms 
a sphere in $\F \times I$, which must be embedded since 
$\interior(D) \cap \X = \emptyset$.
Since $\F \times I$ is irreducible, the sphere bounds a ball, which can 
be used to isotope $\X$ relative to $\de \X$.  If $\de D
\cap \de \X = \emptyset$, then $D'$ can be pushed entirely 
off of $\de D$, and one component of $\X \cap A$ is 
eliminated.  Components of $\X \cap (A \cap A')$, components of 
$\X \cap A'$,  and additional 
components of $\X \cap A$ may also be removed if $\interior(D') 
\cap (A \cup A') \neq \emptyset$, but no new components of any kind 
are added. Therefore, $\zeta(\X)$ goes down. 
If, instead, $\de D \cap \de \X \neq
\emptyset$, then the isotopy of $X$ relative to $\de X$ must leave 
$\de D \cap \de \X$ fixed.  But this isotopy still 
decreases the rank of $H_1(\X \cap A)$ without increasing 
the rank of $H_1(\X \cap A')$ or the number of components of 
$\X \cap (A \cap A')$.  

The explanation is analogous if the roles of 
$A$ and $A'$ are interchanged.  
\smallskip

\begin{figure}[htbp]
\epsfysize=2.0in
\centerline{
\epsfbox{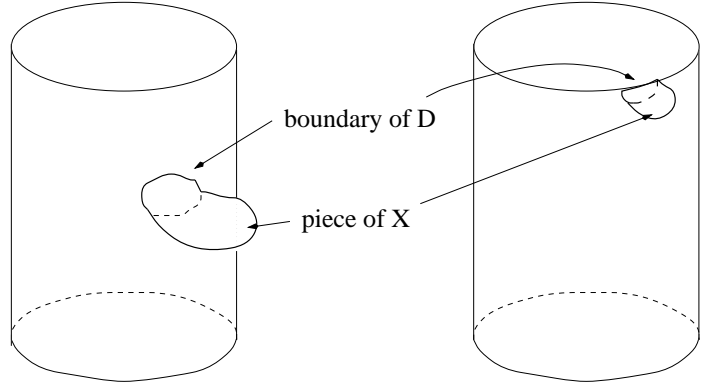}}
\caption{Move 1 disks.} \label{fig-move1}
\end{figure}

\mvmv{Move 2} Suppose E is a disk in $\F \times I$ whose boundary consists
of two arcs $\alpha$ and $\sigma$.  Suppose that $E \cap \X = \sigma$
and that $E \cap A = E \cap A' = \alpha$.  Then $E$ can be used to
isotope $\X$ relative to $\de \X$ and decrease $\zeta(\X)$.  The
isotopy will keep $\X$ in the class of pseudo-transverse embeddings.
\smallskip

\begin{figure}[htbp]
\epsfysize=4.0in
\centerline{
\epsfbox{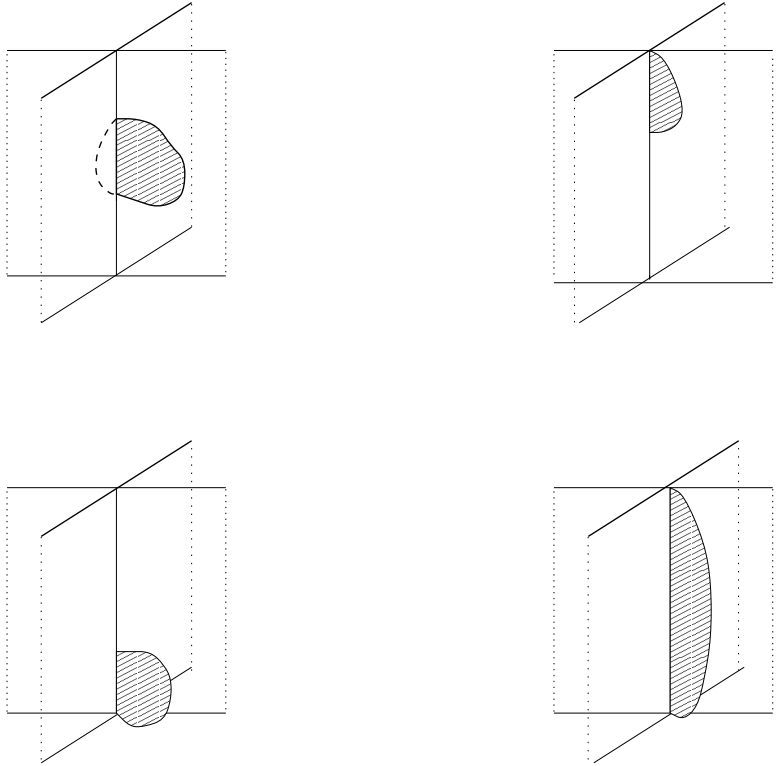}}
\caption{Move 2 disks.} \label{fig-move2}
\end{figure}

\begin{figure}[ptbh]
\epsfysize=7.0in
\centerline{
\epsfbox{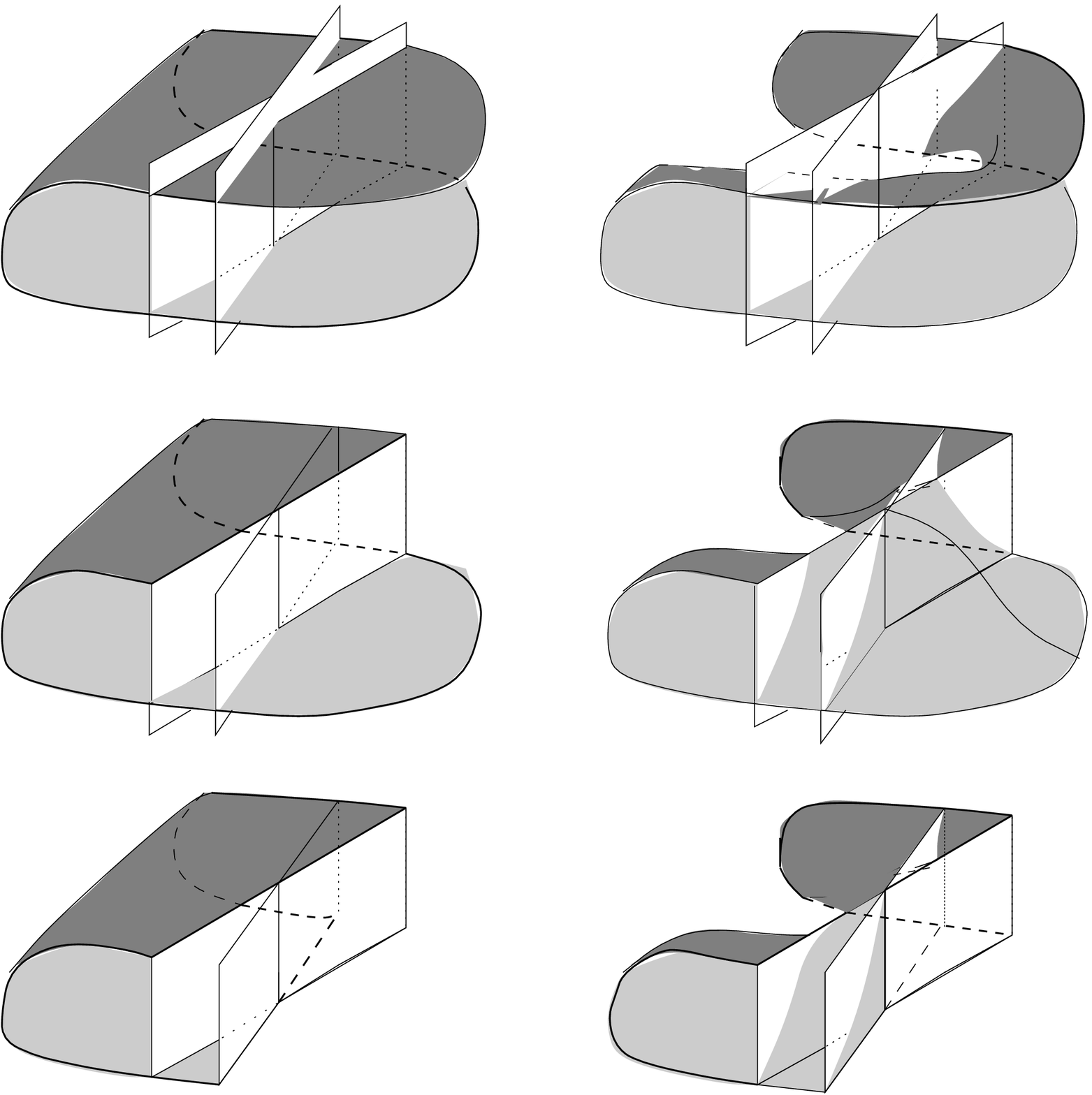}}
\caption{Using move 2 disks to isotope $\X$.} \label{fig-move3}
\end{figure}

\mvmv{Explanation of Move 2}  
Consider three cases depending on how many points of $\de \sigma$ lie in 
$\de \X$.  See Figures~\ref{fig-move2} and~\ref{fig-move3}.

\caca{Case 1} Both endpoints of $\sigma$ lie in $\interior(\X)$.  Then $\X$
can be isotoped in a neighborhood of $\sigma$ so that it moves through
$E$ and slips entirely off of $\alpha$.  This isotopy decreases by two
the number of components of $\X \cap (A \cap A')$.

\caca{Case 2} One endpoint of $\sigma$ lies in $\interior(\X)$ and one endpoint
lies in $\de \X$.  Now $\X$ cannot be isotoped relative to $\de \X$ entirely
off $\alpha$ since the endpoint of $\de \sigma$ in $\de \X$ must
remain fixed.  But $\X$ can still be pushed off $\interior(\alpha)$ and
off the free endpoint, lowering the number of components of $\X
\cap (A \cap A') $ by one.

\caca{Case 3} Both endpoints of $\sigma$ lie in $\de \X$. Notice that one
endpoint must lie in $\X \cap (\F \times 1)$ and one must lie in $\X
\cap (\F \times 0)$, since $\alpha$ is a vertical line connecting them.
In this case $\X$ can be isotoped relative to $\de \X$ in a neighborhood of
$\sigma$ to move $\sigma$ directly onto the vertical line $\alpha$ and
produce a twist around this vertical line resembling picture (3)(i).
This isotopy decreases the number of components of $\X \cap (A \cap
A')$ by one, since it transforms the two endpoints of $\de \sigma$
into a single vertical line.

\mvmv{Move 3} Suppose $W$ is an annulus of $A$ such that 
$W \cap A' = \emptyset$ and $\de W$ is the boundary of an annulus
of $\X \sm A$.  Then $W$ can be used to isotope $\X$ and decrease
$\zeta(\X)$.  The isotopy leaves $\X$ in the class of pseudo-transverse
embeddings.  A similar move is possible if the roles of $A$ and $A'$
are interchanged.

\begin{figure}[htbp]
\epsfysize=2.0in
\centerline{\epsfbox{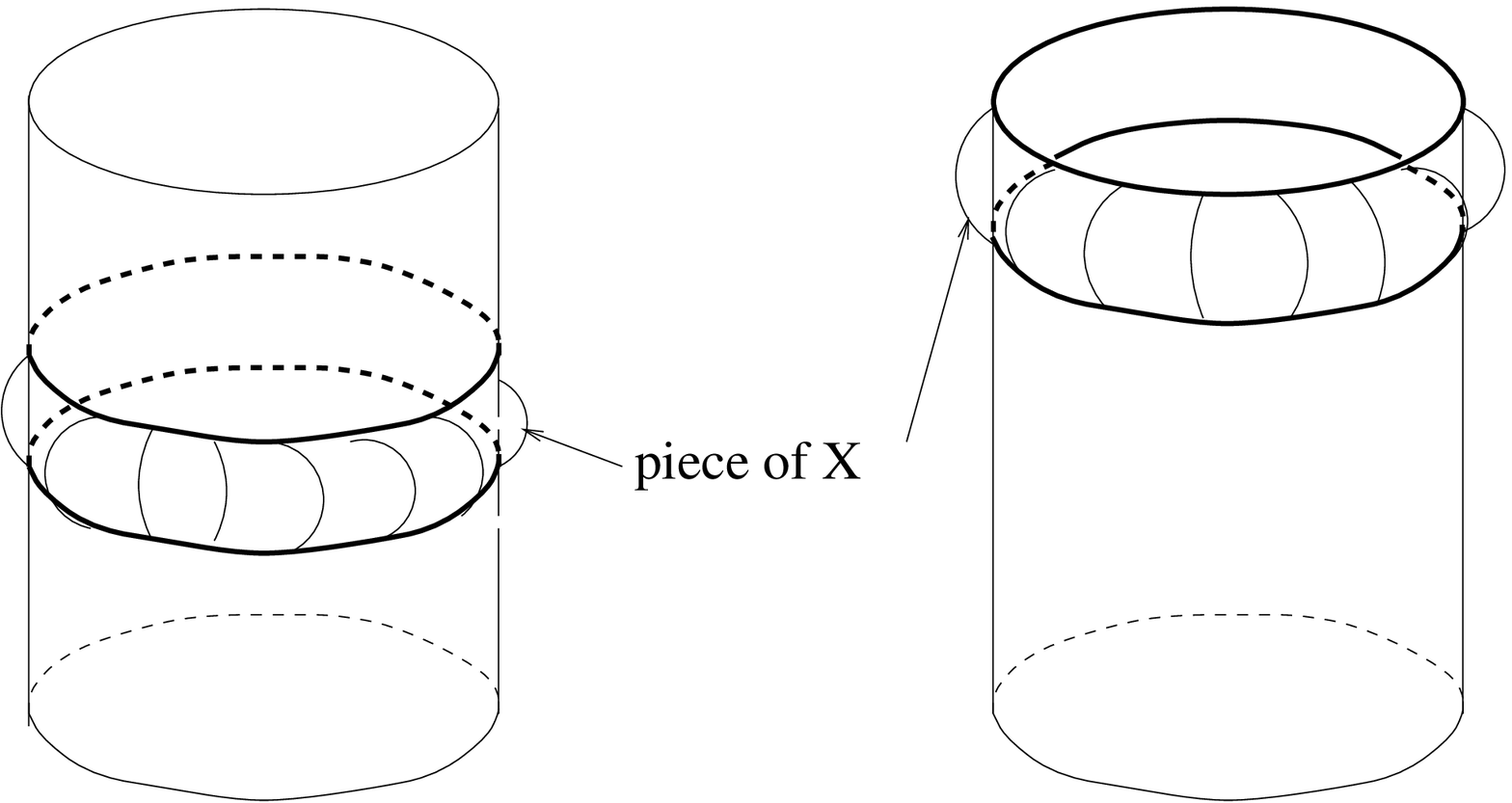}}
\caption{Move 3 annuli.} \label{fig-move3really}
\end{figure}

\mvmv{Explanation of Move 3} Let $K$ be the annulus of $\X \sm A$ such that
$\de K = \de W$, and let $L$ be the closure of the component of $\F
\times I \sm A$ that contains $K$.  $W \cup K$ forms a torus, which is
embedded in $L$ since $\interior(K) \cap A = \emptyset$.  The torus $W
\cup K$ compresses in $L$, because $L$ is homotopic to a surface with
boundary and therefore $\pi_1(L)$ cannot contain a $\Z \times \Z$
subgroup.  Since $L$ is irreducible and $W$ lies on the boundary of
$L$, $W \cup K$ actually bounds a solid torus in $L$, which can be
used to isotope $\X$ relative to $\partial \X$ by pushing $K$ through
$W$. Notice that $\partial W$ and $\de \X$ share at most one
component.  If $\de W$ and  $\de X$ are disjoint, then this isotopy
removes at least two components of $\X \cap A$.  If $\de W$ and $\de
\X$ share a component, the the isotopy removes at least one component of 
$\X \cap A$.  In either case, the isotopy decreases the rank of 
$H_1(\X \cap A)$ without increasing the
rank of $H_1(\X \cap A')$ or the number of components of $\X \cap (A
\cap A')$.

The explanation is analogous if the roles of $A$ and $A'$ are
interchanged.

\medskip
I will refer to the type of disk used in Move 1 as  a
move 1 disk, the type of disk used in Move 2 as a move 2 disk, and the
type of annulus used in Move 3 as a move 3 annulus.

\stst{Step 1} Isotop $\X$ relative to $\de \X$ so that it is
pseudo-transverse to $A \cup A'$.  This can be accomplished, for
example, by making $\X$ honestly transverse to $A \cup A'$.  Then only
pictures (1)(i), (1)(ii), (1)(iv), (2)(i), (2)(ii), (2)(iv), (3)(ii),
(3)(iii), and (3)(iv) in the definition of pseudo-transverse can
occur.

\smallskip
\stst{Step 2} Suppose $\F \times I$ contains a move 1 disk, a move 2
disk, or a move 3 annulus.
Use it to isotope $\X$.  Repeat this step as often as necessary, until
there are no more such disks or annuli. The process must terminate after
finitely many moves, since each move decreases $\zeta(\X)$.  

At this stage,
$\X$ already has a neat posture with respect to $A \cup A'$.  In
particular, the following claims hold, where $K$ is any
component of $\X \sm (A \cup A')$ and $L$ is the component of $(\F
\times I) \sm (A \cup A')$ that contains $K$.

\begin{figure}[htbp]
\epsfxsize=3.0in
\centerline{\epsfbox{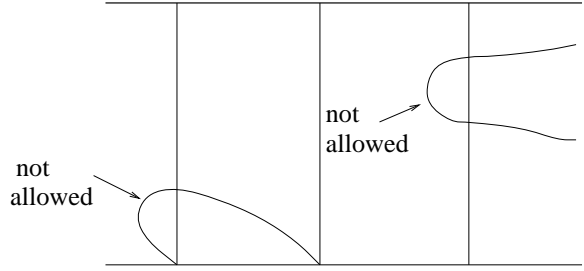}}
\caption{Endpoints of $\mu$ go to distinct vertical lines of $A \cap A'$.} \label{fig-distinct}
\end{figure}

\clcl{Claim 1} Suppose that  $\mu$ is an arc contained in $\X \cap (A
\sm A')$ with both endpoints in $A \cap A'$.  Then either the endpoints of 
$\mu$ go to distinct
vertical lines of $(A \cap A')$ or else $\mu$ wraps all the
way around an annulus of $A$.  Likewise for arcs of $\X \cap (A' \sm
A)$.  See Figure~\ref{fig-distinct}.

\clcl{Claim 2} For any point $z \epsilon (C \sm C')$, there is a neighborhood 
$U$ of $z$ such that $p^{-1}(U) \cap \X$ is a disjoint union of
neighborhoods of the form:
\begin{figure}[h!tbp]
\epsfysize=0.7in
\centerline{\epsfbox{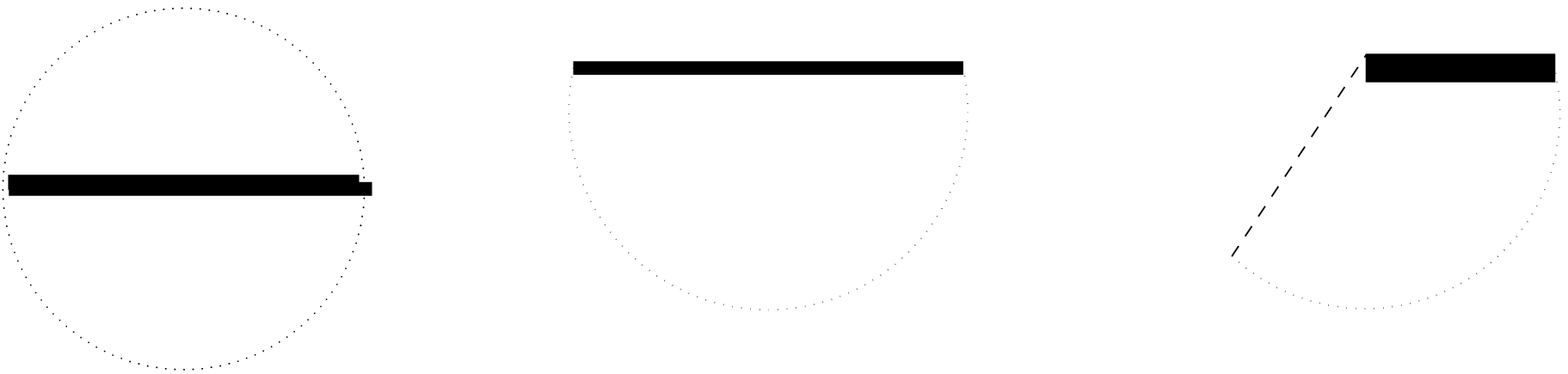}}
\end{figure}

\par\noindent
Likewise for points of $(h(C) \sm C)$. In other words, pictures
(1)(iii), (1)(v), (2)(iii) and (2)(v) in the definition of
pseudo-transverse do not occur.

\clcl{Claim 3}  Every 
circle of $\partial K$ has non-zero degree in the 
cylinder of $\partial_v L$ that contains it.   Here
$\dv L$ refers to the vertical 
boundary of $L$, namely ${\partial L}\cap (A \cup A' \cup (\de \F \times I))$.

\clcl{Claim 4} $\pi_1(K) \rightarrow \pi_1(L)$ is injective.

\clcl{Claim 5} Either $\pi_1(K) \rightarrow \pi_1(L)$ is surjective, or else 
$K$ is an annulus that is parallel to $\dv L$.

\clcl{Claim 6}  If $K$ is an annulus, then the two loops of $\partial K$ 
go to two distinct cylinders of $\dv L$.
\smallskip\noindent

\bigskip

\begin{figure}[htbp]
\epsfxsize=3.0in
\centerline{\epsfbox{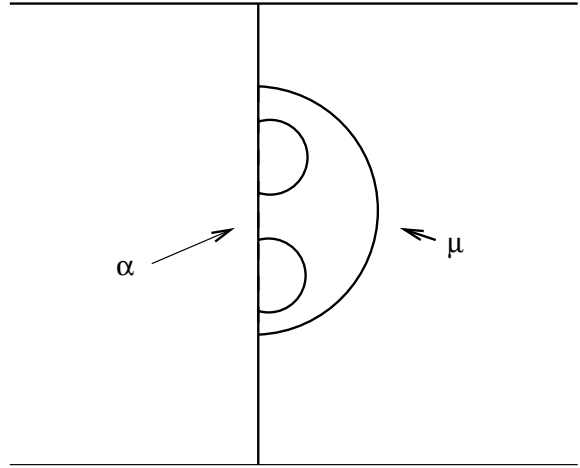}}
\caption{$E \cap \X$ may contain additional arcs.}  \label{fig-arcs}
\end{figure}

\clcl{Proof of Claim 1} Let $\mu$ be an arc contained in $\X \cap (A \sm A')$  and
suppose both endpoints of $\mu$ lie in one vertical line of $A \cap
A'$.  Suppose that $\mu$ does not wrap all the way around an annulus
of $A$, and let $\alpha$ be the segment of $A \cap A'$ that connects
the endpoints of $\mu$.  Then $\alpha \cup \mu$ bounds a half-disk $E$
in $A$.  Notice that $E \cap A' = \alpha$.  The set $E \cap \X$ cannot
contain any closed loops of $\X \cap A$, since any such loop would
bound a move 1 subdisk of E, which should have been removed in Step 2.
But $E \cap \X$ may contain other arcs besides $\mu$ with endpoints on
$\alpha$. (See Figure~\ref{fig-arcs}.)  By replacing $\mu$ and $E$
with an arc and subdisk closer to $\alpha$ if necessary, I can assume
that $E
\cap \X = \mu$.  Nudge $E$ relative $\alpha$ off of $A$ to get a new disk $E'$
bounded by the arcs $\alpha$ and $\mu'$, where $\mu' \subset \X$ and
$\interior(\mu') \subset \interior(\X)$.  Since $E \cap \X =
\mu$ and $E \cap A' = \alpha$,  I can assume that $E' \cap \X =
\mu'$ and $E' \cap A' = \alpha$.  Also, $E' \cap A = \alpha$.  So
$E'$ is a move 2 disk, in violation of Step 2.  Thus, the endpoints of
$\mu$ must lie in distinct components of $A \cap A'$ after all.

The same argument applies to arcs contained in 
$\X \cap (A' \sm A)$.
\medskip

\begin{figure}[htbp]
\epsfxsize=5.0in
\centerline{\epsfbox{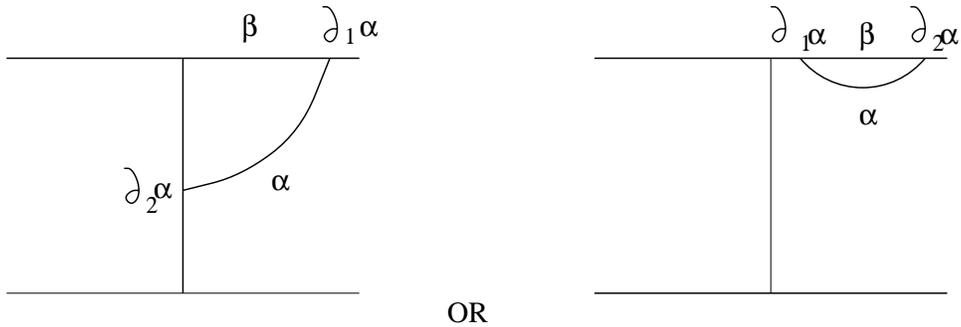}}
\caption{The arc $\alpha$ in the proof of Claim 2.}  \label{fig-ignore}
\end{figure}

\clcl{Proof of Claim 2} Suppose picture (1)(iii) or (1)(v) does occur.
(The argument is similar if picture (2)(iii) or (2)(v) occurs.)  Let
$\alpha$ be the segment of ${{p^{-1}}(U)} \cap \X \cap A$ 
drawn vertically in these pictures,
extended in $A$ until it first hits $\de \X$ or $A \cap A'$.
See Figure~\ref{fig-ignore}.  Label the first endpoint of $\alpha$ as
$\partial_1\alpha$ and the second endpoint as $\partial_2 \alpha$. So
$\partial_1 \alpha$ lies on $C \times 1$.

Suppose first that $\partial_2 \alpha$ lies on $\partial \X$ but not
on $A \cap A'$.  Recall that $\de \X = (C \times 1) \cup (C' \times 0)
\cup B$; therefore $(\de \X \cap A) \sm (A \cap A') \subset ((C \times 1)
\cap A) \cup (B \cap A)$.  Since $p|_B$ is an embedding, it follows
that $B \cap A \subset C \times 1$.  So $\partial_2 \alpha$ lies in $\de
\X \cap (C \times 1))$.  Thus, $\alpha$ cuts off an arc $\beta$ of
$C \times 1$ such that $\alpha \cup \beta$ bounds a move 1 disk.
This disk should already have been removed in Step 2.

Next, suppose that $\partial_2\alpha$ lies on $A \cap A'$, and let
$\beta$ be the arc of $\partial \X$ such that $\partial_1\beta =
\partial_1\alpha$ and $\partial_2\beta$ lies on the same vertical line
of $A \cap A'$ as $\partial_2\alpha$, and so that $\alpha \cup \beta$
does not wrap all the way around an annulus of $A$.  Then $\alpha \cup
\beta$ is an arc contained in $\X \cap A$ with both endpoints in the same
vertical line which should not exist by Claim 1.
\medskip

\begin{figure}[htbp]
\epsfysize=1.7in
\centerline{\epsfbox{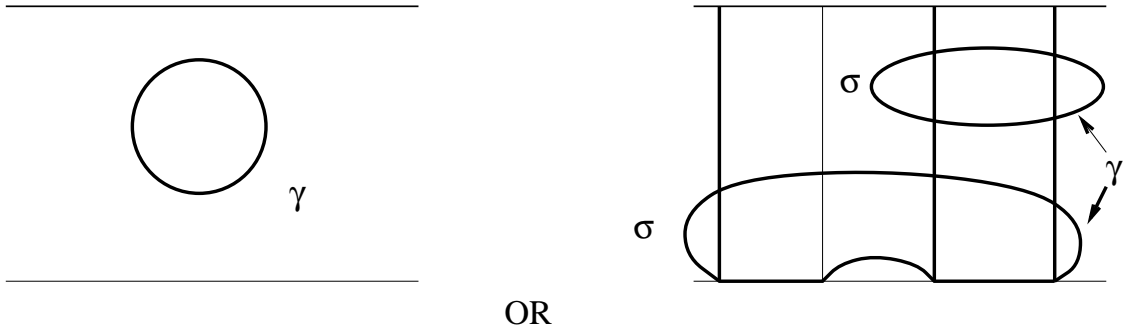}}
\caption{Suppose a circle $\gamma$ of $\partial K$ gets sent to a
cylinder of $\dv L$ with degree 0.}  \label{fig-degree0}
\end{figure}

\clcl{Proof of Claim 3}
Let $K$ be a component of $\X \sm (A \cup A')$ and let $L$ be the 
component of $(\F \times I) \sm (A \cup A')$ that contains it. Let
$\dv L$ denote the vertical boundary of $L$, that is, $\dv L = \de L
\cap (A \cup A' \cup (\de \F \times I))$.  Suppose 
that a circle $\gamma$ of $\partial K$ gets sent to a cylinder $G$ of
$\dv L$ by degree 0.  $G$ may be an annulus of $A$ or an annulus of $A'$.
Or $G$ may consist of rectangles of $A \sm A'$ and $A' \sm
A$, joined together along vertical lines of $A \cap A'$ and possibly
along vertical strips of $\de \F \times I$.  Suppose first that $\gamma$ is
contained in a single component of $A \sm A'$ or $A' \sm A$.  (This
happens, in particular, if $G$ is an annulus of $A$ or $A'$.)  Since
$\gamma$ has degree 0 in $G$, it must bound a disk in $G$, which I can
assume has interior disjoint from $\X$ by replacing $\gamma$ with an
innermost loop if necessary.  But this disk is a move 1 disk, so it
should have been removed already in Step 2.  See
Figure~\ref{fig-degree0}.

Now suppose that $\gamma$ is not contained in a single component of 
$A \sm A'$ or $A' \sm A$.  
Since $p|_B$ is an embedding, $\gamma$ must be disjoint from $\de \F
\times I$.    Again, $\gamma$ bounds a disk  in $G$.
The vertical lines of $A \cap A'$ cut the disk into subdisks.
Consider an outermost subdisk and the arc $\sigma \subset \gamma$ that
forms half its boundary.  This arc $\sigma$ is contained in $\X \cap
(A \sm A')$ (or in $\X \cap (A' \sm A)$) and both its endpoints lie
in the same vertical line.  Furthermore, $\sigma$ cannot wrap all the
way around an annulus of $A$ (or $A'$).  Claim 1 says that such arcs
do not exist.
\medskip

\clcl{Proof of Claim 4}
The following diagram commutes, and the map $\pi_1(\X) \rightarrow
\pi_1(\F \times I)$ is injective.  So it will suffice to show that
$\pi_1(K) \rightarrow \pi_1(\X)$ is injective.

\begin{picture}(150, 80)(-30, 0)
\thinlines
\put(20, 60){$\pi_1(K)$}
\put(80, 60){$\pi_1(L)$}
\put(140, 60){$\pi_1(\F \times I)$}
\put(80, 20){$\pi_1(\X)$}
\put(38, 43){\vector(3,-2){23}}
\put(117, 27){\vector(3,2){25}}
\put(50, 62){\vector(1,0){17}}
\put(110, 62){\vector(1,0){17}}
\end{picture}

If $\pi_1(K) \rightarrow \pi_1(\X)$ does not inject, then some loop
in $\partial K$ must bound a disk $D$ in $\X \sm K$.  I claim
that $D$ contains at least one loop $\omega$ of $\X \cap A$ or $\X
\cap A'$. 

Since $\de D$ bounds $D$ on one side and $K$ on the other, $\de D$
must be disjoint from $\de X$.  So $\de D$ must intersect $A$ or $A'$.
Assume without loss of generality that $X$ intersects $A$.  Possibly
$\partial D$ itself is contained in $\X \cap A$.
In this case set $\omega = \partial D$.  Otherwise, take an arc of
$\partial D \cap (\X \cap A)$  and let
$\omega$ be the component of $\X \cap A$ that
contains it.  Notice that $\omega$ must lie in $D$, since it cannot
intersect $\interior(K)$.  Since $D \cap \de \X = \emptyset$, $\omega$
is a closed loop rather than an arc.

Since $\omega \subset D$, $\omega$ is null-homotopic in $\F \times I$.
Since each annulus in $A$ $\pi_1$-injects into $\F \times I$, $\omega$
must bound a disk $E$ in $A$.  If $E \cap A' = \emptyset$, then $E$ is
a move 1 disk, which should already have been removed.  If $E \cap A'
\neq \emptyset$, then $\de E \cap (A \sm A')$ will contain an arc
whose endpoints lie in the same vertical line of $A \cap A'$, as in
the proof of Claim 3. But Claim 1 says that such arcs do not exist.

\medskip
\clcl{Proof of Claim 5} 
Let $R$ be the surface $p(L)$.  From Claim 4, the map $\pi_1(K) \rightarrow
\pi_1(L)$ is injective.  Therefore, it is possible to write $\pi_1(L)$
either as an amalgamated product or as an HNN extension over
$\pi_1(K)$, depending on whether or not $K$ separates $L$.  
In fact,
$K$ separates $L$, because the inclusion of $R \times 1$
into $L$ induces an isomorphism of fundamental groups, which could not
happen if $\pi_1(L)$ were an HNN extension.  Let $M_1$ be the
component of $L \sm K$ that contains $R \times 1$, and let $M_2$ be the
other component.  Then the composition of maps
\[\pi_1(R \times 1) \rightarrow \pi_1(M_1) \rightarrow\pi_1(L)\] 
is an isomorphism.  So $\pi_1(M_1) \rightarrow \pi_1(L)$ is
surjective.  But \[\pi_1(L) = \pi_1(M_1)\ast_{\pi_1(K)}\pi_1(M_2)\] so
the injection $\pi_1(K) \rightarrow \pi_1(M_2)$ must be an
isomorphism.  By the h-cobordism theorem \cite[Theorem 10.2]{Hempel},
$(M_2, K) \cong (K \times I, K \times 1)$.

Now there are two possibilities: either $R \times 0 \subset M_1$ or $R
\times 0 \subset M_2$.  If $R \times 0 \subset M_1$, then $\partial
M_2 \sm K$ is a subset of $\dv L$.  So $K$
is parallel to a subsurface of $\dv L$, and therefore must be a disk or
an annulus.  By Claim 3, $K$ must be a boundary parallel annulus.

If $(R \times 0) \subset M_2$, then it follows as above that the
composition 
\[\pi_1(R \times 0) \rightarrow \pi_1(M_2) \rightarrow \pi_1(L) \] 
is an isomorphism.  Therefore, $\pi_1(M_2) \rightarrow \pi_1(L)$ is
surjective.  Also, $\pi_1(M_2) \rightarrow \pi_1(L)$ is injective
since 
\[\pi_1(L) = \pi_1(M_1)\ast_{\pi_1(K)}\pi_1(M_2)\]  
So $\pi_1(M_2) \rightarrow \pi_1(L)$ is an isomorphism.  In addition, $\pi_1(K)
\rightarrow \pi_1(M_2)$ is an isomorphism, from above.  
Therefore, $\pi_1(K) \rightarrow \pi_1(L)$ is an isomorphism, and the
claim is proved.

\medskip

\clcl{Proof of Claim 6}
Let $K$ be any annulus component of $\X \sm (A \cup A')$, let $L$ be
the corresponding component of $(\F \times I) \sm (A \cup A')$, and
suppose that both circles of $\partial K$ go to the same cylinder $G$
of $\partial L$.  Since both circles have degree $\pm 1$ in
$\partial_v L$, they bound an annulus $W$ in $\dv L$.  
Notice that $W$ is disjoint from $\de \F \times
I$, since by assumption, $p|_B: B \rightarrow \de \F$ is an embedding.
If $W$ is entirely 
contained in $A$ or $A'$, then $W$ is a move 3 annulus, which should
have been removed in Step 2.  
So $W$ must consist of alternating rectangles of $A$ and $A'$.  

The union $W \cup K$ forms a torus, which is embedded in $L$ since
$\interior(K) \cap \partial_v L = \emptyset$.  The torus $W \cup K$
compresses in $L$, because $L$ is homotopic to a surface with boundary
and therefore $\pi_1(L)$ cannot contain a $\Z \times \Z$ subgroup.
Since $L$ is irreducible and $W$ lies on the boundary of $L$, $W \cup
K$ bounds a solid torus $T$ in $L$.

I will construct a move
2 disk as follows.  Start with a vertical arc $\alpha$ of $W \cap (A \cap
A')$.  Connect its endpoints with an embedded arc $\sigma$ of $K$ so
that $\alpha \cup \sigma$ is null homotopic in $T$. This is possible
because each component of $\partial K$ generates $\pi_1 (T)$, so I can
replace a poor choice of $\sigma$ by one that wraps around $\partial
K$ an additional number of times and get a good choice of $\sigma$.
Now $\alpha \cup \sigma$ bounds an embedded disk $E$ in $T$, which I
can assume has interior disjoint from $\X$ by replacing it with a
subdisk if necessary.  In addition, $E \cap A = E \cap A' = \alpha$,
so $E$ is a move 2 disk. But Step 2 already eliminated all disks of
this form.
\bigskip

\stst{Step 3} Recall that $p: {\F \times I} \rightarrow \F$ is the projection
map.  I will isotope $\X$ relative to $\de \X$ so that for any arc or
loop $\gamma$ of $\X \cap (A \sm A')$ or $\X \cap (A' \sm A)$
$p|_\gamma$ is a local homeomorphism onto its image.  I will consider arcs
first and loops next.

Take any arc $\mu$ of $\X \cap (A \sm A')$.  If $\interior(\mu) \cap
\de \X \neq \emptyset$, then $\mu \subset \de \X$ by Claim 2.  The map
$p$ is already injective on arcs of $\de \X \cap A$ and $\de \X \cap
A'$, so I can leave $\mu$ alone.  If $\interior(\mu) \cap \de \X =
\emptyset$, then isotope $\mu$ relative to $\partial \mu$ to an
embedded arc $\mu'$ in $A \sm A'$ such that $p|_{\mu'}$ is a
homeomorphism onto its image.  This is possible because by Claim 1,
either $\mu$ wraps all the way around an annulus of $A$ or else the
endpoints of $\mu$ lie in distinct vertical lines of $A \cap A'$.
The isotopy can be done in such a way that $\mu'$ does not intersect
any other arcs of $\X \cap (A \sm A')$ that may lie in the same
vertical rectangle.  Perturb $\X$ in a neighborhood of $\mu$ to extend
the isotopy on $\mu$.

Pick another arc of $\X \cap (A \sm A')$ and repeat the 
procedure.  When all the arcs of $\X \cap (A \sm A')$ have 
been pulled taut, continue with arcs of $\X \cap (A' \sm A)$.  

Next, consider any loop $\lambda$ of $\X \cap A$ that does not
intersect $A'$.  By Claim 3, the loop $\lambda$ has degree 1 in $A$,
so it can be isotoped to a loop $\lambda'$ such that $p|_{\lambda}$ an
homeomorphism onto its image.  As before, the isotopy can be done in
such a way that $\lambda'$ does not intersect any other loops of $\X
\cap A$, and the isotopy can be extended to a neighborhood of $\lambda$ in
$\X$.  Loops of $\X \cap A'$ that are disjoint from $A$ can be
isotoped similarly.

I claim that at this stage, for any component $K$ of $X \sm (A \cup
A')$, $p|_{\de K}$ is a local homeomorphism except along vertical lines of
$K \cap (A \cap A')$.  Every point of $\de K$ is either a point on the
interior of an arc of $\de \X \cap (\de \F \times I)$, a point of $\X
\cap A \cap (\de \F \times I)$ or $\X \cap (A' \cap (\de \F \times
I))$, a point of $A \cap A'$, or a point on the interior of an arc or
loop of $\X \cap (A \sm A')$ or $\X \cap (A' \sm A)$.  The map
$p|_{\de K}$ is a local homeomorphism near the first type of point by
assumption.  It is a local homeomorphism near the second type of point
by Claim 2.  It is a local homeomorphism near the third type of point
(away from vertical twist lines) because
$X$ is pseudo-transverse.  Finally, $p|_{\de K}$ is a local
homeomorphism near the fourth type of point by the work done in Step 3.

\smallskip

\stst{Step 4} I will finish isotoping $\X$ so that $p$ is locally
injective except at twist lines.  I will do this by using a 
fact about maps between surfaces.  

\stst{Fact} 
{\it Let $f: (G, \de G) \rightarrow (H, \de H)$ be a map
such that $f|\de G$ is a local homeomorphism and $f_\ast: \pi_1(G)
\rightarrow \pi_1(H)$ is injective.  Then there is a homotopy
$f_\tau:(G, \de G) \rightarrow (H, \de H)$, with $\tau \in I$, $f_0 =
f$, and $f_\tau|_{\de G} = f_0|_{\de G}$ for all $\tau$, such that
either (1) or (2) holds:
\begin{enumerate}   
\item  $G$ is an annulus or Mobius band and $f_1(G) \subset \de H$, or
\item $f_1: G \rightarrow H$ is a covering map
\end{enumerate}
}

The case when $G$ is a disk is easy to verify; all other cases are
covered by \cite[Theorem 13.1]{Hempel}.

Pick a component $K$
of $\X \sm (A \cup A')$, and let $L$ be the corresponding component of
$(\F \times I) \sm (A \cup A')$.  Assume first that $\de K$ does not
contain any vertical arcs of $A \cap A'$.  By Claim 4, the map $\pi_1(K)
\rightarrow \pi_1(L)$ is injective.  Since $p_\ast: \pi_1(L) \rightarrow
\pi_1(p(L))$ is an isomorphism, the composition $(p|_K)_\ast: \pi_1(K)
\rightarrow \pi_1(p(L))$ is injective. Furthermore, by the discussion
following Step 3, $p|{\de K}: \de K \rightarrow \de p(L)$ is a local
homeomorphism on $\de K$.  Therefore, there is a homotopy $f_\tau: K
\rightarrow p(L)$ with $f_0 = p|_K$ and  $f_\tau|_{\de K} = p|_{\de K}$
such that either $K$ is an annulus or Mobius band and $f_1(K) \subset
\de p(L)$, or $f_1$ is a covering map.

If $K$ is
a Mobius band and $f_1(K) \subset \de p(L)$, then $f_0(\de K)$ is a
degree 2 loop in $\de p(L)$ which is impossible since $\de K \subset
\dv L$ is embedded.  By Claim 6, it not possible for $K$ to be an
annulus and $f_1(K)$ a subset of $\de p(L)$.  Therefore, $f_1$
must be a covering map.  By Claims 5 and 6, the map
$\pi_1(K) \rightarrow \pi_1(L)$ is surjective.  So the map
$(p_K)_\ast: \pi_1(K) \rightarrow \pi_1(p(L))$ is surjective.
Therefore, $f_1$ must be a homeomorphism.  

If $\de K$ contains vertical arcs of $A \cap A'$, then $p|_{\de K}$ is
still very close to a local homeomorphism -- in fact, if $\hat{K}$ is
the surface obtained by collapsing each vertical arc of $\de K \cap (A
\cap A')$ to a point, then $p|_{K}$ factors through a map $\hat {p}:
\hat{K} \rightarrow p(L)$ such that $\hat{p}|_{\de \hat{K}}$ is a local
homeomorphism.  So it is still possible to homotope $p|_K$ relative to
$\de K$ to a map $f_1$ such that ${f_1}|_{\interior(K)}$ is a homeomorphism.

The homotopy $f_\tau$ of $K$ relative to 
$\de K$ in $p(L)$
induces a homotopy of $K$ relative to $\de K$ in $L$ which keeps the
vertical coordinate of each point of $K$ constant and changes its
horizontal coordinate according to $f_\tau$.  At the end of the homotopy,
the new surface $K'$ itself will be embedded in $\F \times I$,
since $p|_{K'}$ is a homeomorphism.

Homotope $\X$ as described above for every other component of $\X \sm
(A \cup A')$.  At this stage, each component of $\X \sm (A \cup A')$
is embedded in $\F \times I$.  But two components $K'_1$ and $K'_2$ in
the same piece $L$ of $(\F \times I) \sm (A \cup A')$ might intersect
each other.  If that happens, an additional homotopy of $K$ can be
tacked on to clear up the problem, as follows.

Let $K_1 \subset L$ and $K_2 \subset L$ be two components of $\X \sm (A
\cup A')$ before the homotopy of Step 4 and let $K'_1$ and $K'_2$
be these components after the homotopy.  Suppose that $K'_1$ intersects
the component $K'_2$.  Since $K_1$ and $K_2$ are disjoint, the
component of $L \sm K_1$ that meets $\F \times 0$ either contains
all of $K_2$ or else contains no part of $K_2$.  Therefore, either all
loops of $K_1$ lie above the corresponding loops of $K_2$ or else they
all lie below the corresponding loops.  Since the homotopies of $K_1$
and $K_2$ did not move $\de K_1$ and $\de K_2$, the same statement
holds for loops of $\de K'_1$ and $\de K'_2$.  Therefore, it is
possible to alter the vertical coordinates of $K'_1$ and $K'_2$ to
make the two surfaces parallel, so that $K'_1$ lies entirely above
$K'_2$ or entirely below $K'_2$.  Therefore all
components of $\X \sm (A \cup A')$ can be assumed disjoint after
the homotopy of Step 4.

In its final position, $\X$ is embedded in $\F \times I$.  A theorem
of Waldhausen \cite[Corollary 5.5]{W} states that if $G$ and $H$ are
incompressible surfaces embedded in an irreducible 3-manifold, and
there is a homotopy from $G$ to $H$ that fixes $\de G$, then there is
an isotopy from $G$ to $H$ that fixes $\de G$.  Therefore, the above
sequence of homotopies can be replaced by an isotopy.

\section{Remarks}

It is not true that every surface in near-horizontal position is
incompressible.  For example, the surface in Figure~\ref{fig-example}
compresses in $\R^2 \times I$.  I hope to find
simple conditions for when a surface in
near-horizontal position is incompressible.  A special case is considered 
in \cite{thesis}.

\bibliography{bib}

\bibliographystyle{plain}

\end{document}